\input amstex\documentstyle{amsppt}  
\pagewidth{12.5cm}\pageheight{19cm}\magnification\magstep1
\topmatter
\title{Positive conjugacy classes in Weyl groups}\endtitle
\author G. Lusztig\endauthor
\address{Department of Mathematics, M.I.T., Cambridge, MA 02139}\endaddress
\thanks{Supported by NSF grant DMS-1566618.}\endthanks
\endtopmatter   
\document

\define\da{\dagger}

\define\Irr{\text{\rm Irr}}

\define\mpb{\medpagebreak}

\define\op{\oplus}
   
\define\part{\partial}
\define\emp{\emptyset}

\define\n{\notin}

\define\m{\mapsto}
\define\do{\dots}

\define\bsl{\backslash}

\define\sub{\subset}    

\define\T{\times}
\define\ti{\tilde}
\define\nl{\newline}
\redefine\i{^{-1}}

\define\ov{\overline}
\define\ot{\otimes}
\define\bbq{\bar{\QQ}_l}

\define\tr{\text{\rm tr}}

\define\ph{\phi}

\define\r{\rho}

\redefine\l{\lambda}

\define\x{\xi}

\define\Ph{\Phi}

\define\kk{\bold k}

\define\qq{\bold q}

\define\NN{\bold N}

\define\QQ{\bold Q}

\define\ZZ{\bold Z}

\define\cd{\Cal D}

\define\cl{\Cal L}

\define\cn{\Cal N}

\define\car{\Cal R}

\define\ta{\ti a}

\define\sha{\sharp}

\define\BOU{Bo}
\define\DL{DL}
\define\FT{FT}
\define\GP{GP}
\define\CHEV{Ch}
\define\COX{L1}
\define\REP{L2}
\define\ORA{L3}
\define\SP{Sp}

\subhead 1\endsubhead
Let $W$ be a Weyl group. In this paper we introduce the notion of positive 
conjugacy class of $W$. This generalizes the notion of elliptic regular conjugacy class 
in the sense of Springer \cite{\SP}.

Let $w\m|w|$ be the length function on $W$. Let $S=\{s\in W;|s|=1\}$.

Let $v$ be an indeterminate. Recall that the Iwahori-Hecke algebra of $W$ is
the associative $\QQ(v)$-algebra $H$ which, as a $\QQ(v)$-vector space has
basis $\{T_w;w\in W\}$ and has multiplication given by
$T_wT_{w'}=T_{ww'}$ if $|ww'|=|w|+|w'|$ and
$(T_s+1)(T_s-v^2)=0$ if $s\in S$; note that $T_1$ is the unit element of $H$.
This is a split semisimple algebra. Let $\qq=v^2$.

For $w,w'$ in $W$ let $N^{w,w'}$ be the trace of the
$\QQ(v)$-linear map $H@>>>H$, $h\m T_whT_{w'{}\i}$.   
We have $N^{w,w'}\in\ZZ[\qq]$ and
$$N^{w,w'}=\sum_{E\in\Irr W}\tr(T_w,E_v)\tr(T_{w'},E_v)\tag a$$
where $\Irr W$ is the set of irreducible $\QQ[W]$-modules up to isomorphism
and for $E\in\Irr W$, $E_v$ denotes the corresponding simple $H$-module.
Note that when $v$ is specialized to $1$, $H$ becomes the group algebra
$\QQ[W]$ of $W$ and $N^{w,w'}$ specializes to
$N^{w,w'}(1)$, the number of elements $y\in W$ such that $wy=yw'$.
In particular, if $w\in W$, $N^{w,w}$ specializes to
$n^w$, the order of the centralizer of $w$ in $W$; thus the
polynomial $N^{w,w}$ can be viewed as a $\qq$-analogue of the number $n^w$.

If $C$ is a conjugacy class in $W$ we denote by $C_{min}$ the set of all
$y\in C$ such that $C@>>>\NN$, $w\m|w|$, reaches its minimum at $y$. By a result of Geck
and Pfeiffer \cite{\GP,3.2.9}, for any $E\in\Irr W$, $w\m\tr(T_w,E_v)$ is constant on 
$C_{min}$. Using this 
and (a) we see that $w\m N^{w,w}$ is constant on $C_{min}$. We say that
$C$ is {\it positive} if $C\ne\{1\}$ and for some/any $w\in C_{min}$ we have 
$N^{w,w}\in\NN[\qq]$. (We then also say that any element $w\in C_{min}$ is positive.)

\subhead 2\endsubhead
For any $f\in\ZZ[\qq]$ we write $f=\sum_{i\ge0}f_i\qq^i$ where $f_i\in\ZZ$.
For $w\in W$ let $S_w$ be the set of all $s\in S$ such that $s$ appears in some/any
reduced expression for $w$. Let $\cl_w=\{s\in S;|sw|<|w|\}$, $\car_w=\{s\in S;|ws|<|w|\}$.
For $a,a'$ in $W$ we write $T_aT_{a'}=\sum_{b\in W}\ph(a,a',b)T_b$
where $\ph(a,a',b)\in\ZZ[\qq]$. We show:

(a) {\it $\ph(a,a',b)_{|a|}\ge0$, $\ph(a,a',b)_i=0$ for $i>|a|$. If 
$\ph(a,a',b)_{|a|}\ne0$ then either $a'=b$, $S_a\sub\cl_{a'}$ or $|b|<|a'|$.
If $a'=b$, $S_a\sub\cl_{a'}$ then $\ph(a,a',b)_{|a|}=1$.}
\nl
We argue by induction on $|a|$. When $|a|=0$ the result is obvious. Assume now that
$|a|\ge1$. We write $a=a_1s$ where $s\in S$, $|a_1|=|a|-1$. If $|sa'|=|a'|+1$ then
$\ph(a,a',b)=\ph(a_1,sa',b)$ and the induction hypothesis shows that
$\ph(a,a',b)_i=0$ for $i\ge|a|$. Since $s\in S_a$, $s\n\cl_{a'}$, we see that the desired
result holds. Next we assume that $|sa'|=|a'|-1$. Then
$\ph(a,a',b)=\qq\ph(a_1,sa',b)+(\qq-1)\ph(a_1,a',b)$.
From the induction hypothesis we see that $\ph(a,a',b)_i=0$ if $i>|a|$, that 
$\ph(a,a',b)_{|a|}=\ph(a_1,sa',b)_{|a_1|}+\ph(a_1,a',b)_{|a_1|}\ge0$
and that if $\ph(a,a',b)_{|a|}\ne0$ then either $\ph(a_1,sa',b)_{|a_1|}\ne0$
or $\ph(a_1,a',b)_{|a_1|}\ne0$, so that we are in one of the cases (i)-(iv) below.

(i) $sa'=b$,

(ii) $|b|<|sa'|$;

(iii) $a'=b$, $S_{a_1}\sub\cl_{a'}$;

(iv) $|b|<|a'|$.

In case (i),(ii),(iii) we have $|b|<|a'|$; in case (iii) have $a'=b$ and 
$S_a\sub\cl_{a'}$ (since $S_a=S_{a_1}\cup\{s\}$), as desired.
Now assume that
$a'=b$, $S_a\sub\cl_{a'}$. 
It remains to show that $\ph(a_1,sa',b)_{|a_1|}+\ph(a_1,a',b)_{|a_1|}=1$.
By the induction hypothesis we have
$\ph(a_1,a',b)_{|a_1|}=1$ (since $S_{a_1}\sub \cl_{a'}$)
$\ph(a_1,sa',b)_{|a_1|}=0$ (since $sa'\ne b$ and $|b|\not<|sa'|$).
This completes the proof.
\nl
The following result is proved in the same way as (a).

(b) {\it $\ph(a,a',b)_{|a'|}\ge0$, $\ph(a,a',b)_i=0$ for $i>|a'|$. If 
$\ph(a,a',b)_{|a'|}\ne0$ then either $a=b$, $S_{a'}\sub\car_{a'}$ or $|b|<|a|$.
If $a=b$, $S_{a'}\sub\car_{a'}$ then $\ph(a,a',b)_{|a'|}=1$.}
\nl
For $a,a',a''$ in $W$ we have $T_aT_{a'}T_{a''}=\sum_{b\in W}f(a,a',a'',b)T_b$
where $f(a,a',a'',b)\in\ZZ[\qq]$. Let $n=|a|+|a''|$. We show:

(c) {\it $f(a,a',a'',a')_n\ge0$, $f(a,a',a'',a')_i=0$ for $i>n$. If 
$f(a,a',a'',a')_n\ne0$ then $S_a\sub\cl_{a'}$ and $S_{a''}\sub\car_{a'}$. Conversely, if
$S_a\sub\cl_{a'}$ and $S_{a''}\sub\car_{a'}$ then $f(a,a',a'',a')_n=1$.}
\nl
We have $f(a,a',a'',a')=\sum_{c\in W}\ph(a,a',c)\ph(c,a'',a')$. Hence for $i\ge0$ we have
$f(a,a',a'',a')_i=\sum_{c\in W;j\ge0,j'\ge0,j+j'=i}\ph(a,a',c)_j\ph(c,a'',a')_{j'}$. 
Using (a) and (b) in the last sum we can take $j\ge|a|,j'\ge|a''|$. Hence if 
$i>n=|a|+|a''|$ then $f(a,a',a'',a')_i=0$ and 
$f(a,a',a'',a')_n=\sum_{c\in W}\ph(a,a',c)_{|a|}\ph(c,a'',a')_{|a''|}\ge0$.
Assume now that $f(a,a',a'',a')_n\ne0$. Then in the last sum we can assume that

$c=a'$, $S_a\sub\cl_{a'}$ or $|c|<|a'|$ and
$a'=c$, $S_{a''}\sub\car_c$ or $|a'|<|c|$
 
Thus we can assume that $c=a'$, $S_a\sub\cl_{a'}$ and $S_{a''}\sub\car_{a'}$ and using
again (a),(b) we have $f(a,a',a'',a')_n=1$. 

\mpb

For $w,w'$ in $W$ we set $n=|w|+|w'|$; we show:

(d) {\it $N^{w,w'}_n=\sha(a'\in W;S_w\sub\cl_{a'}, S_{w'}\sub\car_{a'})>0$,
$N^{w,w'}_i=0$ for $i>n$.}
\nl
We have $N^{w,w'}=\sum_{a'\in W}f(w,a',w',a')$ and the result follows from (c). (We use 
that if $a'=w_0$, the longest element of $W$ then $S_w\sub\cl_{a'}=S$, 
$S_{w'}\sub\car_{a'}=S$.)

\mpb

From (d) we deduce:

(e) {\it Let $w,w',n$ be as in (d). Assume that either $S_w=S$ or $S_{w'}=S$. then
$N^{w,w'}_n=1$.}
\nl
Indeed, if $a'\in W$ satisfies $S\sub\cl_{a'}$ or $S\sub\car_{a'}$ then $a'=w_0$.

\mpb

We state the following result.

(f) {\it Let $w,w',n$ be as in (d). For $i=0,1,\do,n$ we have
$N^{w,w'}_i=(-1)^nN^{w,w'}_{n-i}$.}
\nl
Let $\bar{}:\QQ(v)@>>>\QQ(v)$ be the field automorphism such that $\bar{v}=v\i$.
For $E\in\Irr W$ let $E^\da\in\Irr W$ be the tensor
product of $E$ with the sign representation of $W$. It is known that for $w\in W$ we have
$$\tr(T_w,E_v^\da)=(-v^2)^{|w|}\tr(T_{w\i}\i,E_v)=(-v^2)^{|w|}\ov{\tr(T_w,E_v)}.$$ 
It follows that
$$\align&N^{w,w'}=\sum_{E\in\Irr W}\tr(T_w,E_v^\da)\tr(T_{w'},E_v^\da)\\&=
\sum_{E\in\Irr W}(-v^2)^{|w|}\ov{\tr(T_w,E_v)}(-v^2)^{|w'|}\ov{\tr(T_{w'},E_v)}\\&
=\sum_{E\in\Irr W}(-v^2)^n\ov{\tr(T_w,E_v)}\ov{\tr(T_{w'},E_v)}.\endalign$$
We see that 
$$N^{w,w'}=(-\qq)^n\ov{N^{w,w'}}$$
and (f) follows.

\subhead 3\endsubhead
We now assume that $W$ is irreducible. 
Let $\nu=|w_0|$ where $w_0$ is the longest element of $W$.
An element $w\in W$ (or its conjugacy class) is said to be {\it elliptic} if its
eigenvalues in the reflection representation of $W$ are all $\ne1$. 
For any $d\in\{2,3,4,\do\}$ let $C^d$ be the set
of all elliptic elements $w\in W$ which have order $d$ and are regular in the sense of
Springer \cite{\SP}. It is known \cite{\SP} that $C^d$ is either empty or a single 
conjugacy class in $W$. Let $\cd=\{d\in\{2,3,\do\}; C^d\ne\emp\}$. It is known \cite{\SP}
that if $d\in\cd$ and $w\in C^d_{min}$ then $|w|=2\nu/d$. Let $h$ be the Coxeter number of
 $W$. We have $h\in\cd$.

According to \cite{\SP}, the set $\cd$ is as follows: 

Type $A_n(n\ge1)$: $\cd=\{n+1\}$.

Type $B_n (n\ge2)$: $\cd=\{d\in\{2,4,6,\do\};2n/d=\text{integer}\}$.

Type $D_n$ ($n$ even, $n\ge4$): 

$\cd=\{d\in\{2,4,6,\do\};(2n-2)/d=\text{odd integer or} 2n/d=\text{integer}\}$.

Type $D_n$ ($n$ odd, $n\ge5$): $\cd=\{d\in\{2,4,6,\do\};(2n-2)/d=\text{odd integer}\}$.

Type $E_6$: $\cd=\{3,6,9,12\}$.

Type $E_7$: $\cd=\{2,6,14,18\}$.

Type $E_8$: $\cd=\{2,3,4,5,6,8,10,12,15,20,24,30\}$.

Type $F_4$: $\cd=\{2,3,4,6,8,12\}$.

Type $G_2$: $\cd=\{2,3,6\}$.
\nl
We note the following properties:

(a) If $2\in\cd$, $d\in\cd$ is even and $w\in C^d_{min}$ then $w^{d/2}=w_0$, 
$(d/2)|w|=|w_0|$ hence $T_w^{d/2}=T_{w_0}$.

(b) If $d=h$, $w\in C^d_{min}$ then $T_w^d=T_{w_0}^2$. 

(c) If $d\in\cd$, $h/d\in\NN$ and $y\in C^h_{min}$ then $y^{h/d}\in C^d_{min}$ and
$(h/2)|y|=|y^{h/d}|$ hence $T_y^{h/d}=T_{y^{h/d}}$.
\nl
The equation $w^{d/2}=w_0$ in (a) holds by examining the characteristic polynomial
of $w$ and $w^{d/2}$ in the reflection representation of $W$; then (a) follows.
The equality in (b) can be deduced from \cite{\BOU, Ch.V,\S6, Ex.2}.
The equation $w^{h/d}\in C^d_{min}$ in (c) holds by examining the characteristic polynomial
of $w$ and $w^{h/d}$ in the reflection representation of $W$; then (c) follows.

\mpb

For any $E\in \Irr W$ we define $a_E\in\NN$ as in \cite{\ORA, 4.1}.
Let $\ta_E=\nu-a_E+a_{E^\da}$. 

(d) $T_{w_0}^2=v^{2\ta_E}1:E_v@>>>E_v$.
\nl
This can be deduced from \cite{\ORA, (8.12.2)}; a closely related statement was first
proved by Springer, see \cite{\GP, 9.2.2}.

We show:

(e) {\it Let $E\in\Irr W$ and let $d\in\cd,w\in C^d_{min}$. Then all eigenvalues 
of $T_w:E_v@>>>E_v$ (in an algebraic closure of $\QQ(v)$) are roots of $1$ times
$v^{2\ta_E/d}$.}
\nl
If $d$ is as in (a) then the result follows from (a) and (d). If $d=h$ then the result 
follows from (b) and (d). 
If $d,y$ are as in (c) then the result follows from (c) and the previous sentence. 
From the description of $\cd$ for various types we see that if $d\in\cd$ is not as in (a)
then it is as in (c). This proves (e). (A closely related result can be found in
\cite{\GP, 9.2.5}.)

\mpb

From (e) we deduce:

(f) {\it In the setup of (e), $\tr(T_w,E_v)$ equals $v^{2\ta_E/d}\tr(w,E)$; this is $0$ 
if $2\ta_E/d\n\ZZ$.}
\nl
(The idea of the proof leading to (f) appeared in \cite{\ORA, p.320}.) Using (f) and 1(a) 
we deduce:

(g) {\it If $d\in\cd$, $w\in C^d_{min}$, then 
$$N^{w,w}=\sum_{E\in\Irr W}\qq^{2\ta_E/d}\tr(w,E)^2.$$
In particular, we have $N^{w,w}\in\NN[\qq]$ and $w$ is positive.}

\subhead 4\endsubhead
Using 1(a) and the CHEVIE package \cite{\CHEV} one can find a list of positive 
conjugacy classes in $W$ (assumed to be irreducible of low rank). I thank Gongqin Li for 
help with programming in GAP. The list of positive conjugacy classes in $W$
which are not regular elliptic for $W$ of type 
$E_6,E_7,E_8,F_4,G_2,B_5,B_6$ is as follows. (We specify a
conjugacy class by the characteristic polynomial of one of its elements in the reflection
representation. We denote by $\Ph_k$ the $k$-th cyclotomic polynomial; thus $\Ph_2=\qq+1$,
$\Ph_3=\qq^2+\qq+1$, etc.)

Type $E_6$: none.

Type $E_7$: $\Ph_{12}\Ph_6\Ph_2,\Ph_{10}\Ph_6\Ph_2,\Ph_{10}\Ph_2^3,\Ph_8\Ph_4\Ph_2$,
$\Ph_4^2\Ph_2^3$.

Type $E_8$: $\Ph_{18}\Ph_6,\Ph_{18}\Ph_2^2,\Ph_9\Ph_3,\Ph_{14}\Ph_2^2$.

Type $F_4$: none.

Type $G_2$: none.

Type $B_5$: $\Ph_8\Ph_2,\Ph_6\Ph_2^2,\Ph_4^2\Ph_2,\Ph_2\Ph_4\Ph_6$.

Type $B_6$: $\Ph_{10}\Ph_2^2,\Ph_8\Ph_4,\Ph_8\Ph_2^2,\Ph_6\Ph_2^3$.
\nl
In each of these examples any positive element of $W$ is elliptic; we expect this to 
be true in general. The example of $B_6$ suggests that if $W$ is of type $B_n$ with 
$2n=4+8+\do+4k$, then an 
element of $W$ with cycle type $(4)(8)\do(4k)$ might be positive.

{\it Remark.}
In a first version of this paper, the fourth conjugacy class listed
above for type $B_5$ was omitted by mistake. I thank Jean Michel for pointing this
out.

\subhead 5\endsubhead
Let $\kk$ be an algebraic closure of the finite field $F_q$ with $q$ elements. 
Let $G$ be a connected reductive group over $\kk$ with a fixed
$F_q$-split rational structure and whose Weyl group is $W$. Let $F:G@>>>G$ be the
corresponding Frobenius map. For $w\in W$ let $X_w$ be the variety of Borel subgroups
$B$ of $G$ such that $B$ and $F(B)$ are in relative position $w$, see \cite{\DL, 1.3}. The 
finite group $G^F=\{g\in G;F(g)=g\}$ acts on $X_w$ by conjugation. For $w,w'$ in $W$
we denote by $X_{w,w'}=G^F\bsl(X_w\T X_{w'})$ the space of $G^F$-orbits for the diagonal 
action of $G^F$
on $X_w\T X_{w'}$. Now $(B,B')\m(F(B),F(B'))$ induces a map $X_{w,w'}@>>>X_{w,w'}$ 
(denoted again by $F$) which is the Frobenius map for an $F_q$-rational structure on 
$X_{w,w'}$. By \cite{\REP, 3.8} for any integer $e\ge1$ we have
$$\sha(\x\in X_{w,w'};F^e(\x)=\x)=N^{w,w'}(q^e).\tag a$$

\subhead 6\endsubhead
In the remainder of this paper we assume that $G$ in no.5 is simply connected
and $W$ is irreducible.
In the case where $w$ is a Coxeter element of minimal length of $W$, the left hand side
of 5(a) (with $w=w'$) has been computed in \cite{\COX, p.158}. This gives the following
formulas for $N^{w,w}$.

Type $A_n (n\ge1)$: $\qq^{2n}+\qq^{2n-2}+\do+\qq^2+1$.

Type $B_n (n\ge2)$: $\qq^{2n}+2\qq^{2n-2}+2\qq^{2n-4}+\do+2\qq^2+1$.

Type $D_n (n\ge4)$: $\qq^{2n}+\qq^{2n-2}+2\qq^{2n-4}+2\qq^{n-6}+\do+2\qq^4+\qq^2+1$.

Type $E_6$: $\qq^{12}+\qq^{10}+2\qq^8+4\qq^6+2\qq^4+\qq^2+1$.

Type $E_7$: $\qq^{14}+\qq^{12}+2\qq^{10}+4\qq^8+2\qq^7+4\qq^6+2\qq^4+\qq^2+1$.

Type $E_8$: $\qq^{16}+\qq^{14}+2\qq^{12}+4\qq^{10}+2\qq^9+10\qq^8+2\qq^7+4\qq^6+2\qq^4
+\qq^2+1$.

Type $F_4$: $\qq^8+2\qq^6+6\qq^4+2\qq^2+1$.

Type $G_2$: $\qq^4+4\qq^2+1$.
\nl
Let $\cn_G$ be the variety consisting of all pairs
$(g,g')$ where $g$ runs through the standard Steinberg cross section of the set of regular
elements of $G$ and $g'$ is an element in the centralizer of $g$ in $G$ modulo the centre 
of $G$. (This variety, introduced in \cite{\COX, p.158}, makes sense even if $\kk$ is 
replaced by the complex numbers. It plays a role in \cite{FT} where it is called the
{\it universal centralizer}.)
According to \cite{\COX, p.158}, the number of $F_q$-rational points of $\cn_G$
is equal to $N^{w,w}(q)$ hence it is given by the formulas above with $\qq=q$.

\subhead 7\endsubhead
Let $C$ be a conjugacy class of $W$.
For $w\in C$, the part of weight $j$ of the $i$-th $l$-adic cohomology 
space with compact support $H^i_c(X_w,\bbq)$
is a direct sum $\op_{\r}V^i_{\r,j}\ot\r$ where 
$\r$ runs over the unipotent representations of $G^F$ (up to isomorphism) and
$V^i_{\r,\j}$ are finite dimensional $\bbq$-vector spaces in such a way that the 
$G^F$-action is only through the action on $\r$ and the Frobenius action is only through
an action on $V^i_{\r,j}$ (where it is multiplication by $q^{j/2}\l_\r$ with
$\l_\r$ a root of $1$ independent of $w,i,j$, and the parity of $j$ is independent
of $w,i$, see \cite{\REP, 3.9}, \cite{\ORA}). Using
the Grothendieck-Lefschetz fixed point formula, from 5(a) we deduce for any $e\ge1$:
$$N^{w,w}(q^e)=\sum_{i,i',j,j',\r}(-1)^{i+i'}\dim(V^i_{\r,j})\dim(V^{i'}_{\r^*,j'})
q^{je/2}q^{j'e/2}$$
where $\r^*$ is the dual of $\r$ and we have used that $\l_{\r^*}=\l_\r\i$. This implies
$$N^{w,w}=\sum_{i,i',j,j',\r}(-1)^{i+i'}\dim(V^i_{\r,j})\dim(V^{i'}_{\r^*,j'})v^{j+j'}.\tag a$$
If we assume that 

(b) the $G^F$-modules $H^i_c(X_w,\bbq)$, $\text{ dual of }H^{i'}_c(X_w,\bbq)$ are
disjoint for any $i,i'$ such that $i\ne i'\mod2$
\nl
then from (a) we could deduce that $N^{w,w}\in\NN[\qq]$. Hence if we assume further
that $w\in C_{min},C\ne\{1\}$ it would follow that $C$ is positive.

We conjecture that, conversely, if $C$ is positive and $w\in C_{min}$ then
(b) holds. It is also likely that in this case,

(c) the $G^F$-modules $H^i_c(X_w,\bbq)$, $H^{i'}_c(X_w,\bbq)$ are
disjoint for any $i,i'$ such that $i\ne i'\mod2$.
\nl
This disjointness property holds when $w$ is as in \S6, see
\cite{\COX}.

\enddocument